\documentclass[fleqn,12pt]{article}
\usepackage{amsfonts}
\usepackage{amsthm}
\usepackage{amssymb}
\usepackage{amsmath}

\newtheorem{theo}{\bf Theorem}[section]
\newtheorem{lemma}{\bf Lemma}[section]
\newtheorem{coro}{\bf Corollary}[section]
\newtheorem{rem}{\bf Remark}[section]

\newcommand{\Q}{{\Bbb Q}}
\newcommand{\R}{{\mathbb R}}

\newcommand{\bea}{\begin{eqnarray*}}
\newcommand{\eea}{\end{eqnarray*}}
\newcommand{\be}{\begin{eqnarray}}
\newcommand{\ee}{\end{eqnarray}}

\newcommand{\vol}{\mbox{vol}\,}

\newcommand{\ve}{\bf}

\newcommand{\GCD}{\mbox{\rm GCD}\,}

\begin{document}

\setlength{\unitlength}{1mm}
\setcounter{page}{1}

\title{An Optimal Lower Bound for the Frobenius Problem}
\author{by Iskander Aliev (Edinburgh) and Peter M. Gruber (Vienna)}
\date{}
\maketitle

\noindent {\bf Abstract.} Given $N\ge 2$ positive integers $a_1,
a_2, \ldots, a_N$ with $\GCD(a_1, \ldots, a_N)=1$, let
$f_N$ denote the largest natural number which is not a positive
integer combination of $a_1, \ldots, a_N$. This paper gives an
optimal lower bound for $f_N$ in terms of the absolute inhomogeneous
minimum of the standard $(N-1)$-simplex.


\bigskip
\noindent {\bf 2000 MS Classification:} 11D85, 11H31, 52C17

\section{Introduction and statement of results}

Given $N\ge 2$ positive integers $a_1, a_2, \ldots, a_N$ with
$\GCD(a_1, \ldots, a_N)=1$, the Frobenius problem asks for the
largest natural number $g_N=g_N(a_1, \ldots, a_N)$ (called the
Frobenius number) such that $g_N$ has no representation as a
non--negative integer combination of $a_1, \ldots, a_N$. In this
paper, without loss of generality,  we assume that
$a_1<a_2<\ldots<a_N$. The simple statement of the Frobenius problem
makes it attractive and the relevant bibliography is very large (see
\cite{Alf} and Problem C7 in \cite {Guy}). We will mention just few
main results.
%

For $N=2$, the Frobenius number is given by an explicit formula due
to W. J. Curran Sharp \cite{Sharp}: 
\bea g_2(a_1,a_2)=(a_1-1)(a_2-1)-1\,. \eea
The case $N=3$ was solved explicitly by Selmer and Beyer \cite{SB},
using a continued fraction algorithm. Their result was simplified by
R\"odseth \cite{Ro} and later by Greenberg~\cite{Gre}.
No general formulas are known for $N\ge 4$.  Upper bounds, among
many others, include classical results by Erd\"os and Graham
\cite{EG}
\bea g_N\le 2 a_N\left [\frac{a_1}{N}\right ]-a_1\,, \eea
by Selmer \cite{Se}
\bea g_N\le 2 a_{N-1}\left [\frac{a_N}{N}\right ]-a_N\,, \eea
and by Vitek \cite{Vi}
\bea g_N\le \left [\frac{(a_2-1)(a_N-2)}{2}\right ]-1\,, \eea
as well as more recent results by Beck, Diaz, and Robins \cite{BDR}
\bea g_N\le \frac{1}{2}\left(\sqrt{a_1 a_2
a_3(a_1+a_2+a_3)}-a_1-a_2-a_3\right)\,, \eea
and by Fukshansky and Robins \cite{Lenny}, who produced an upper
bound
in terms of the covering radius of a lattice related to the integers
$a_1, \ldots, a_N$.

For $N=3$, Davison \cite{Dav} has found a sharp lower bound
\bea g_3\ge \sqrt{3a_1 a_2 a_3}-a_1-a_2-a_3\,, \eea
where the constant $\sqrt{3}$ cannot be replaced by any smaller
constant.
R\"odseth \cite{Ro} proved in the general case that
\bea g_N\ge ((N-1)! a_1 \cdots a_{N} )^{1/(N-1)}-\sum_{i=1}^N a_i\,.
\eea
%

The present paper gives a sharp lower bound for the function
\bea f_N(a_1, \ldots, a_N)=g_N(a_1, \ldots, a_N)+ \sum_{i=1}^N
a_i\eea
(and thus for $g_N$) in terms of a geometric characteristics of the
standard $(N-1)$-simplex. Clearly, $f_N=f_N(a_1, \ldots, a_N)$ is
the largest integer which is not a {\em positive} integer
combination of $a_1, \ldots, a_N$.

Following the geometric approach developed in \cite{Kannan} and
\cite{Kannan-Lovasz}, we will make use of tools from the geometry of
numbers.
Recall that a
family of sets in
$\R^{N-1}$ is a {\em covering} if their union equals $\R^{N-1}$.
Given a 
set $S$ and a lattice $L$, we say that $L$ is a {\em covering
lattice} for
$S$ 
if the family $\{S+{\ve l}: {\ve l}\in L\}$
%
is a covering. 
%
%
Recall also that the {\em inhomogeneous minimum} of 
the set $S$ with respect to
the lattice $L$ is the quantity
\bea \mu(S,L)= \inf\{\sigma>0: L \,\,\text{a covering lattice
of}\,\, \sigma S\} \eea
and the quantity
\bea \mu_0(S)=\inf\{\mu(S,L): \det L =1\}\, \eea
is called the {\em absolute inhomogeneous minimum} of $S$. If $S$ is
bounded and has inner points,
then $\mu_0(S)$ does not vanish and is finite (see \cite{GrLek},
Chapter 3).

Let $S_{N-1}$ be the standard simplex given by
\bea S_{N-1}=\{(x_1,\ldots,x_{N-1}): x_i\ge 0 \,\mbox{reals and}\,
\sum_{i=1}^{N-1}x_i \le 1\}\,.\eea
The main result of the paper shows that the constant
$\mu_0(S_{N-1})$ is a sharp lower bound for (suitably normalized)
Frobenius number and integers with relatively small $f_N$ are,
roughly speaking, dense in $\R^{N-1}$.
\begin{theo}
%
%
\begin{itemize}
\item[{\rm (i)}] For $N\ge 3$ the inequality
\be \mu_0(S_{N-1})\le \frac{f_N(a_1, \ldots, a_N)}{(a_1 \cdots
a_{N})^{1/(N-1)}}\,\label{main_inequality}\ee
holds.
\item[{\rm (ii)}]
For any $\epsilon>0$ and for any point ${\ve
\alpha}=(\alpha_1,\ldots,\alpha_{N-1})$ in $\R^{N-1}$ there exist
$N$ integers $0<a_1<a_2< \ldots< a_N$ with $\GCD(a_1, \ldots,
a_N)=1$ such that
\be \left | \,\alpha_i-\frac{a_i}{a_N}\right |<\epsilon
\,,\;\;\;i=1,2,\ldots,N-1\label{Density}\ee
and
\be \frac{f_N(a_1, \ldots, a_N)}{(a_1 \cdots a_{N})^{1/(N-1)}}<
\mu_0(S_{N-1})+\epsilon\,. \label{Sharpness}\ee
\end{itemize}
\label{Th_of_Kannan}
\end{theo}
\begin{rem}
{\rm Prof. L. Davison kindly informed the authors that the part (i)
of Theorem \ref{Th_of_Kannan} was proved by R\"odseth in \cite{Ro1}
without using geometry of numbers.}
\end{rem}

%
%
%

%
The quantity $\mu_0(S)$ is closely related to the {\em covering
constant} $\Gamma(S)$ of the set $S$, where
\be \Gamma(S)= \sup\{\det(L): L \,\,\text{a covering lattice of}\,\,
S\}\,. \label{Gamma}\ee
By Theorem 1, Ch. 3, \S 21 of \cite{GrLek} (see also \cite{Bambah})
for each Lebesgue measurable set $S$
\be \Gamma(S)\le \vol(S)\,, \label{gamma_vol}\ee
and by Theorem 2 ibid.
\be \mu_0(S)=\frac{1}{\Gamma(S)^{1/(N-1)}}\,. \label{gamma_mu}\ee
The proof of Theorem 1, Ch. 3, \S 21 of \cite{GrLek} easily implies
that the equality in (\ref{gamma_vol}) is attended only if $S$ is a
space--filler. Further, by Theorem 6.3 of \cite{Rogers}, packings of
simplices cannot be very dense and, consequently, $S_{N-1}$ is not a
space--filler.
Therefore, by (\ref{gamma_vol}) and (\ref{gamma_mu}),
\be \mu_0(S_{N-1})> \frac{1}{(\vol
(S_{N-1}))^{1/(N-1)}}=((N-1)!)^{1/(N-1)}\,, \label{Lower_S}\ee
and we get the following result.
\begin{coro} For $N\ge 3$ the inequality
\be f_N(a_1, \ldots, a_N)> ((N-1)! a_1 \cdots a_{N} )^{1/(N-1)}\,
\label{Kb}\ee
holds. \label{Kill_coro}
\end{coro}
The inequality (\ref{Kb}) with non--strict sign was proved in
\cite{Ro1}.
%
%
%
%
%
%
%
%
%
The only known value of $\mu_0(S_{N-1})$ is $\mu_0(S_2)=\sqrt{3}$
(see e. g. \cite{Fary}). In the latter case we get the following
slight generalization of Theorems 2.2 and 2.3 in \cite{Dav}.
\begin{coro} For $N=3$ the inequality
\bea f_3(a_1, a_2, a_3)\ge  (3 a_1 a_2 a_3)^{1/2}\eea
holds. Moreover, for any $\epsilon>0$ and for any point ${\ve
\alpha}=(\alpha_1,\alpha_2)$ in $\R^{2}$ there exist integers
$0<a_1< a_2< a_3$ with $\GCD(a_1, a_2, a_3)=1$ such that
\bea \left | \,\alpha_i-\frac{a_i}{a_3}\right |<\epsilon
\,,\;\;\;i=1,2\eea
and
\bea f_3(a_1, a_2, a_3)< ((3+\epsilon)a_1 a_2 a_3)^{1/2}\,. \eea

\label{Result_of_Davison}
\end{coro}
Let us consider a lattice $M$ in $\R^{N-1}$ generated by the vectors
\be \frac{1}{N-1} {\ve e}_1, \ldots,  \frac{1}{N-1} {\ve
e}_{N-1}\,,\label{covering_lattice}\ee
where ${\ve e}_j$ are the standard basis vectors. Since the
fundamental cell of $M$ w. r. t. the basis (\ref{covering_lattice})
belongs to $S_{N-1}$, the lattice $M$ is a covering lattice for the
simplex $S_{N-1}$. Therefore, by (\ref{Gamma})
and (\ref{gamma_mu}), 
\bea \mu_0(S_{N-1})\le \frac{1}{(\det M)^{1/(N-1)}}=N-1\,.
\eea
This implies the following result.
\begin{coro}
For any $\epsilon>0$ and for any point ${\ve
\alpha}=(\alpha_1,\ldots,\alpha_{N-1})$ in $\R^{N-1}$ there exist
$N$ integers $0<a_1<a_2< \ldots< a_N$ with $\GCD(a_1, \ldots,
a_N)=1$ such that
\bea \left | \,\alpha_i-\frac{a_i}{a_N}\right |<\epsilon
\,,\;\;\;i=1,2,\ldots,N-1\eea
and
\bea \frac{f_N(a_1, \ldots, a_N)}{(a_1 \cdots a_{N})^{1/(N-1)}}<
N-1+\epsilon\,. \eea
\label{N-1}
\end{coro}

\begin{rem} {\rm Note that the inequality (\ref{Lower_S}) and Stirling's
formula imply that}
\bea \liminf_{N\rightarrow\infty}\frac{\mu_0(S_{N-1})}{N-1}\ge
e^{-1}\,. \eea
{\rm Thus, 
we know the asymptotic behavior of
the optimal constant $\mu_0(S_{N-1})$ up to the multiple $e$}.
\end{rem}
For ${\ve a}=(a_1, a_2, \ldots , a_N)$, define a lattice $L_{\ve a}$
by
\bea L_{\ve a}=\{(x_1,\ldots,x_{N-1}): x_i \,\mbox{integers and}\,
\sum_{i=1}^{N-1}a_i x_i \equiv 0 \mod a_N\}\,. \eea
The following theorem is implicit in~\cite{NewSL}.
\begin{theo}
For any lattice $L$ with basis ${\ve b}_1,\ldots,{\ve b}_{N-1}$,
${\ve b}_i\in\mathbb Q^{N-1}$, $i=1,\ldots,N-1$ and for all
rationals $\alpha_1, \ldots,\alpha_{N-1}$ with
$0<\alpha_1\le\alpha_2\le\cdots\le\alpha_{N-1}\le 1$, there exists
an infinite arithmetic progression ${\mathcal P}$ and a sequence
\bea {\ve a}(t)=(a_1(t),\ldots,a_{N-1}(t),a_{N}(t))\in\mathbb
Z^{N}\,, t\in\mathcal P\,,\eea such that
$\GCD(a_1(t),\ldots,a_{N-1}(t),a_{N}(t))=1$ and the lattice $L_{{\ve
a}(t)}$ has a basis \bea{\ve b}_1(t),\ldots,{\ve b}_{N-1}(t)\eea
with
\be \frac{b_{ij}(t)}{d\,t}=b_{ij}+O\left(\frac{1}{t}\right)\,,\;\;\;
i,j=1,\ldots,{N-1}\,, \label{asympt_aij} \ee
where $d\in\mathbb N$ is such that $d\, b_{ij}, d\,
\alpha_j\,b_{ij}\in\mathbb Z$ for all $i,j=1,\ldots,{N-1}$.
Moreover,
\be a_{N}(t)=\det(L)d^{N-1}t^{N-1}+O(t^{N-2})\label{asympt_h} \ee
and
\be \alpha_i(t):=\frac{a_i(t)}{a_{N}(t)}
=\alpha_i+O\left(\frac{1}{t}\right). \label{asympt_alpha}\ee
\label{main_lemma}
\end{theo}
For completeness, we give a proof of Theorem \ref{main_lemma} in
Section \ref{Proof_of_the_main_lemma}.

\section{Proof of Theorem \ref{Th_of_Kannan} (i)}
%
%
%
Recall that ${\ve a}=(a_1, a_2, \ldots , a_N)$ and put
\bea
\alpha_1=\frac{a_1}{a_N}\,,\ldots\,,\alpha_{N-1}=\frac{a_{N-1}}{a_N}\,.\eea
Define a simplex $S_{\ve a}$ by
\bea S_{\ve a}=\{(x_1,\ldots,x_{N-1}): x_i\ge 0 \,\mbox{reals and}\,
\sum_{i=1}^{N-1}a_i x_i \le 1\} \,.\eea
Theorem 2.5 of \cite{Kannan} states that
%
%
\be f_N(a_1, \ldots, a_N)=\mu(S_{\ve a}, L_{\ve a})\,.
\label{result_of_Kannan}\ee
Observe that the inhomogeneous minimum $\mu(S, L)$ satisfies
\bea \mu(S, tL)=t\mu(S, L)\,, \eea \bea \mu(tS, L)=t^{-1}\mu(S,
L)\,. \eea
Thus, if we define
\bea S_{\ve \alpha}=a_N S_{\ve a}=\{(x_1,\ldots,x_{N-1}): x_i\ge 0
\,\mbox{reals and}\, \sum_{i=1}^{N-1}\alpha_i x_i \le 1\} \,, \eea
\bea L_u=a_N^{-1/(N-1)}L_{\ve a} \eea
then
%
\be \mu (S_{\ve a}, L_{\ve a})=a_N^{1+1/(N-1)}\mu (S_{\ve \alpha},
L_u)\,.\label{Vynos} \ee
Note that $\det L_{\ve a}=a_N$. Thus the lattice $L_u$ has
determinant $1$ and we have
\be \mu_0(S_{\ve \alpha})\le \mu (S_{\ve \alpha}, L_u)\,.
\label{S_alpha}\ee
The simplices $(\alpha_1 \cdots \alpha_{N-1})^{1/(N-1)}S_{\ve
\alpha}$ and $S_{N-1}$ are equivalent up to a linear transformation
of determinant $1$.
Therefore
\be \mu_0(S_{N-1})= \frac{\mu_0(S_{\ve \alpha})}{(\alpha_1 \cdots
\alpha_{N-1})^{1/(N-1)}}\,, \label{transform}\ee
and by (\ref{S_alpha}), (\ref{Vynos}) and (\ref{result_of_Kannan})
we have
\bea \mu_0(S_{N-1})\le \frac{\mu (S_{\ve \alpha}, L_u)}{(\alpha_1
\cdots \alpha_{N-1})^{1/(N-1)}}=\frac{\mu (S_{\ve a}, L_{\ve
a})}{a_N^{1+1/(N-1)}(\alpha_1 \cdots \alpha_{N-1})^{1/(N-1)}} \eea
\bea = \frac{f_N(a_1,\ldots,a_N)}{(a_1\cdots a_N)^{1/(N-1)}}\,.\eea

\section{Proof of Theorem \ref{Th_of_Kannan} (ii)}

The proof is based on Theorem \ref{main_lemma} and the following
continuity property of the inhomogeneous minima. We say that a
sequence $S_t$ of {\em star bodies} in $\R^{N-1}$ converges to a
star body $S$ if the sequence of {\em distance functions} of $S_t$
converges uniformly on the unit ball in $\R^{N-1}$ to the distance
function of $S$.
\begin{lemma}
Let $S_t$ be a sequence of star bodies in $\mathbb R^{N-1}$ which
converges to a bounded star body $S$ and let $L_t$ be a sequence of
lattices in $\mathbb R^{N-1}$ convergent to a lattice $L$. Then
\bea \lim_{t\rightarrow\infty}\mu(S_t,L_t)=\mu(S,L)\,. \eea
\label{limit}
\end{lemma}
\begin{proof}
The result follows from a much more general Satz 1 of \cite{Peter}.
\end{proof}

W. l. o. g., we may assume that ${\ve \alpha}\in \Q^{N-1}$ and
\be 0<\alpha_1<\alpha_2<\ldots <\alpha_{N-1}<
1\,.\label{conditions_on_alpha} \ee
%
%
%
%
%
%
For $\epsilon> 0$ we can choose a lattice $L_\epsilon$ of
determinant $1$ with
\be \mu(S_{\ve \alpha},L_\epsilon)<\mu_0(S_{\ve
\alpha})+\frac{\epsilon(\alpha_1\cdots\alpha_{N-1})^{1/(N-1)}}{2}\,.\label{Step_1}\ee
%
%
The inhomogeneous minimum is independent of translation and rational
lattices are dense in the space of all lattices. Thus, by Lemma
\ref{limit}, we may assume that $L_\epsilon\subset\Q^{N-1}$.
%
%
%
%
%
%
Applying Theorem \ref{main_lemma} to the lattice $L_\epsilon$ and
the numbers $\alpha_1, \ldots, \alpha_{N-1}$, we get a sequence
${\ve
a}(t)$, 
satisfying (\ref{asympt_aij}), (\ref{asympt_h}) and
(\ref{asympt_alpha}). Note also that, by
(\ref{conditions_on_alpha}),
\bea 0<a_1(t)<a_2(t)<\ldots <a_{N}(t)\,\eea
for sufficiently large $t$.
%

Observe that the inequality (\ref{asympt_alpha}) implies
(\ref{Density}) with $a_i=a_i(t)$, $i=1,\ldots,N$, for $t$ large
enough. Let us show that, for sufficiently large $t$, the inequality
(\ref{Sharpness}) also holds.
%
%
Define a simplex $S_{{\ve \alpha}(t)}$ and a lattice $L_t$ by
%
%
%
\bea S_{{\ve \alpha}(t)}=a_N(t)S_{{\ve
a}(t)}=\{(x_1,\ldots,x_{N-1}): x_i\ge 0 \,\mbox{reals and}\,
\sum_{i=1}^{N-1}\alpha_i(t) x_i \le 1\} \,, \eea
\bea L_t=a_N(t)^{-1/(N-1)}L_{{\ve a}(t)}\,. \eea
%

By (\ref{asympt_aij}) and (\ref{asympt_h}),
the sequence 
$L_t$ converges to the lattice $L_\epsilon$.
%
Next, the point ${\ve p}=(1/(2N), \ldots, 1/(2N))$ is an inner point
of the simplex $S_\alpha$ and all the simplicies $S_{{\ve
\alpha}(t)}$ for sufficiently large $t$. By (\ref{asympt_alpha}) and
Lemma \ref{limit}, the sequence $\mu(S_{{\ve \alpha}(t)}-{\ve p},
L_t)$ converges to $\mu(S_\alpha-{\ve p}, L_\epsilon)$. Since the
inhomogeneous minimum is independent of translation, the sequence
$\mu(S_{{\ve \alpha}(t)}, L_t)$ converges to $\mu(S_\alpha,
L_\epsilon)$. Consequently, by (\ref{asympt_alpha}),
\bea \frac{\mu(S_{{\ve \alpha}(t)},
L_t)}{(\alpha_1(t)\cdots\alpha_{N-1}(t))^{1/(N-1)}} \rightarrow
\frac{\mu(S_\alpha,
L_\epsilon)}{(\alpha_1\cdots\alpha_{N-1})^{1/(N-1)}}\,,\;\;\mbox{as}\;t\rightarrow\infty\,,\eea
and, by (\ref{result_of_Kannan}), 
(\ref{Step_1}) and (\ref{transform}),
\bea \frac{f_N(a_1(t),\ldots ,a_{N}(t))}{(a_1(t)\cdots
a_{N}(t))^{1/(N-1)}}=\frac{\mu (S_{{\ve \alpha}(t)},
L_t)}{(\alpha_1(t)\cdots
\alpha_{N-1}(t))^{1/(N-1)}}<\mu_0(S_{N-1})+\epsilon\,\eea
for sufficiently large $t$.

%
%
%
%
%
%
%
%
%
%
%

\section{Proof of Theorem \ref{main_lemma}}
\label{Proof_of_the_main_lemma}

Let us consider the matrices
\bea B=\left (
\begin{array}{ccccc}
b_{11} & b_{12}  & \ldots & b_{1\,N-1} &
\sum_{i=1}^{N-1}\alpha_ib_{1i}\\
b_{21} & b_{22}  & \ldots & b_{2\,N-1} &
\sum_{i=1}^{N-1}\alpha_ib_{2i}\\
\vdots & \vdots & & \vdots & \vdots\\
b_{N-1\,1} & b_{N-1\,2}  & \ldots & b_{N-1\,N-1} &
\sum_{i=1}^{N-1}\alpha_ib_{N-1\,i}\\
\end{array}
\right ) \eea
and
\bea M=M(t,t_1,\ldots,t_{N-1})\eea\bea =\left (
\begin{array}{ccccc}
db_{11}t+t_1 & db_{12}t  & \ldots & db_{1\,N-1}t &
d\sum_{i=1}^{N-1}\alpha_ib_{1i}t\\
db_{21}t & db_{22}t+t_2  & \ldots & db_{2\,N-1}t &
d\sum_{i=1}^{N-1}\alpha_ib_{2i}t\\
\vdots & \vdots & & \vdots & \vdots\\
db_{N-1\,1}t & db_{N-1\,2}t  & \ldots & db_{N-1\,N-1}t+t_{N-1} &
d\sum_{i=1}^{N-1}\alpha_ib_{N-1\,i}t\\
\end{array}
\right ) \,.\eea
Denote by $M_i=M_i(t,t_1,\ldots,t_{N-1})$ and $B_i$ the minors
obtained by omitting the $i$th column in $M$ or in $B$,
respectively.
Following the proof of Theorem 2 in \cite{NewSL}, we observe that
\be |B_{N}|=|\det(b_{ij})|=\det L\,,\label{last_coordinate} \ee
\be |B_i|=\alpha_i|B_{N}|\,,\label{coordinates} \ee
\be M_i=d^{N-1}B_it^{N-1}+\mbox{polynomial of degree less than}\,\,
N-1\,\, \mbox{in}\,\, t\,, \label{asympt_minors} \ee
and $M_1,\ldots,M_{N}$ have no non--constant common factor.

By Theorem 1 of \cite{NewSL} applied with $m=1$, $F=1$, and
$F_{1\nu}=M_\nu(t,t_1,\ldots,t_{N-1})$, $\nu=1,\ldots,N$,
there exist integers $t^*_1,\ldots,t^*_{N-1}$ and an infinite
arithmetic progression ${\mathcal P}$ such that for $t\in{\mathcal
P}$
\bea\mbox{GCD}(M_1(t,t^*_1,\ldots,t^*_{N-1}), \ldots
,M_{N}(t,t^*_1,\ldots,t^*_{N-1}))=1 \,.\eea
Put
\bea {\ve a}(t)=(M_1(t,t^*_1,\ldots,t^*_{N-1}), \ldots
,(-1)^{N-1}M_{N}(t,t^*_1,\ldots,t^*_{N-1}))\,,\;\;\;t\in{\mathcal
P}\,.\eea
Then the basis ${\ve b}_1(t),\ldots,{\ve b}_{N-1}(t)$ for $L_{{\ve
a}(t)}$ satisfying the statement of Theorem \ref{main_lemma} is
given by the rows of the matrix obtained by omitting the $N$th
column in the matrix $M(t,t^*_1,\ldots,t^*_{N-1})$. The properties
(\ref{last_coordinate})--(\ref{asympt_minors}) of minors $M_i$,
$B_i$ imply the properties (\ref{asympt_aij})--(\ref{asympt_alpha})
of the sequence ${\ve a}(t)$, $t\in{\mathcal P}$.

\vskip .5cm \noindent{\bf Acknowledgement}. The authors are
especially grateful to Professors M. Henk and A. Schinzel for
important comments and remarks that strongly improve the exposition.
The authors also wish to thank Professors I. Cheltsov, L.
Fukshansky, L. Davison and J. Ramirez Alfonsin for very helpful and
useful discussions.

\noindent School of Mathematics, University of Edinburgh,  Edinburgh
EH9 3JZ, UK, I.Aliev@ed.ac.uk

\bigskip
\noindent Technische Universit\"at Wien, Wiedner Hauptstraße 8-10 /
1046, 1040 Wien, Austria, peter.gruber@tuwien.ac.at


\begin{thebibliography}{99}

\bibitem{Bambah}
%
R. P. Bambah, {\em On Lattice Coverings}, Proc. Nat. Inst. Sci.
India, {\bf 19} (1953), 447--459.
%
\bibitem{BDR}
M. Beck, R. Diaz, S. Robins, {\em The Frobenius Problem, Rational
Polytopes, and Fourier-Dedekind sums}, J. Number Theory, {\bf 96}
(2002), no. 1, 1--21.

\bibitem{Sharp}
W. J. Curran Sharp, {\em Solution to Problem 7382} (Mathematics),
Educational Times, {\bf 41} (1884).

\bibitem{Dav}
J.~L.~Davison, {\em On the Linear Diophantine Problem of Frobenius},
J. Number Theory, {\bf 48} (1994), no. 3, 353--363.
%
\bibitem{EG}
P. Erd\"os, R. Graham, {\em On a Linear Diophantine Problem of
Frobenius}, Acta Arith., {\bf 21} (1972), 399--408.


\bibitem{Fary}
I.~F\'{a}ry, {\em Sur la Densit\'{e} des R\'{e}seaux de Domaines
Convexes}, Bull. Soc. Math. France, {\bf 78}(1950), 152--161.
%
\bibitem{Lenny}
L. Fukshansky, S. Robins, {\em Frobenius Problem and the Covering
Radius of a Lattice}, submitted.


\bibitem{Gre}
H. Greenberg, {\em Solution to a Linear Diophantine Equation for
Nonnegative Integers}, J. Algorithms, {\bf 9} (1988), no. 3,
343--353.

%
\bibitem{Guy}
R. K. Guy, Unsolved Problems in Number Theory, Third edition.
Problem Books in Mathematics. Unsolved Problems in Intuitive
Mathematics, Springer-Verlag, New York, 2004.
%

\bibitem{Peter}
P. Gruber, {\em Zur Gitter\"uberdeckung des $\R\sp{n}$ durch
Sternk\"orper}, \"Osterreich. Akad. Wiss. Math.-Natur. Kl. S.-B. II,
{\bf 176} (1967), 1--7.
%
\bibitem{GrLek}
P.~M.~Gruber, C.~G.~Lekkerkerker, {\em Geometry of Numbers},
North--Holland, Amsterdam 1987.
%
%
\bibitem{Kannan}
R. Kannan, {\em Lattice Translates of a Polytope and the Frobenius
Problem}, Combinatorica, {\bf 12}(2)(1992), 161--177.
%
\bibitem{Kannan-Lovasz}
R. Kannan, L. Lov\'asz, {\em Covering Minima and Lattice-Point-Free
Convex Bodies}, Ann. of Math. (2) {\bf 128} (1988), no. 3, 577--602.
%
%
%
\bibitem{Alf}
J. L. Ram\'{\i}rez Alfons\'{\i}n, {\em The Diophantine Frobenius
Problem}, Oxford Lecture Series in Mathematics and Its Applications,
2005.
%
\bibitem{Ro}
O. R\"odseth, {\em On a Linear Diophantine Problem of Frobenius}, J.
Reine Angew. Math., {\bf 301} (1978), 171--178.

\bibitem{Ro1}

O. R\"odseth, {\em An upper bound for the $h$-range of the postage
stamp problem}, Acta Arith., {\bf 54} (1990), no. 4, 301--306.


\bibitem{Rogers}
C. A. Rogers, {\em Packing and Covering}, Cambridge Tracts in
Mathematics and Mathematical Physics, No. 54 Cambridge University
Press, New York 1964.
%
\bibitem{NewSL}
A.~Schinzel, {\em A Property of Polynomials with an Application to
Siegel's Lemma}, Monatsh. Math., {\bf 137} (2002), 239--251.
%
\bibitem{Se}
E. Selmer, {\em  On the Linear Diophantine Problem of Frobenius}, J.
Reine Angew. Math., {\bf 293/294} (1977), 1--17.

\bibitem{SB}
E. Selmer, O. Beyer, {\em On the Linear Diophantine Problem of
Frobenius in Three Variables}, J. Reine Angew. Math., {\bf 301}
(1978), 161--170.


%
\bibitem{Vi}
Y. Vitek, {\em Bounds for a Linear Diophantine Problem of
Frobenius}, J. London Math. Soc., (2) {\bf 10} (1975), 79--85.
%
%
\end{thebibliography}
\end{document}